\documentclass[10pt]{article}%
\usepackage[tbtags]{amsmath}
\usepackage{epsfig,amstext,amssymb,amsthm,latexsym}
\usepackage{amssymb,color}
\usepackage{amsfonts}
\usepackage{amssymb}
\usepackage{graphicx}%
\providecommand{\U}[1]{\protect\rule{.1in}{.1in}}
\providecommand{\U}[1]{\protect\rule{.1in}{.1in}}
\providecommand{\U}[1]{\protect\rule{.1in}{.1in}}
\providecommand{\U}[1]{\protect\rule{.1in}{.1in}}
\providecommand{\U}[1]{\protect\rule{.1in}{.1in}}
\allowdisplaybreaks[4]
\definecolor{c20}{rgb}{0.,0.7,0.}
\definecolor{c30}{rgb}{0.,0.,1.}
\definecolor{c40}{rgb}{1,0.1,0.7}
\definecolor{c50}{rgb}{1,0,0}

\def\rH#1{\textcolor{c20}{#1}}

\def\bE#1{\textcolor{c50}{#1}}
\def\bE#1{#1}

\def\rH#1{#1}

\newtheorem{theo}{Theorem}[section]
\newtheorem{sat}[theo]{Proposition}
\newtheorem{de}[theo]{Definition}
\newtheorem{lem}[theo]{Lemma}
\newtheorem{korr}[theo]{Corollary}

\newtheorem{remarks}[theo]{Remarks}

\newcommand{\BQN}{\begin{eqnarray}}
\newcommand{\EQN}{\end{eqnarray}}
\newcommand{\BQNY}{\begin{eqnarray*}}
\newcommand{\EQNY}{\end{eqnarray*}}
\newcommand{\BS}{\begin{sat}}
\newcommand{\ES}{\end{sat}}
\newcommand{\BT}{\begin{theo}}
\newcommand{\ET}{\end{theo}}
\newcommand{\BK}{\begin{korr}}
\newcommand{\EK}{\end{korr}}

\newcommand{\BD}{\begin{de}}
\newcommand{\ED}{\end{de}}
\newcommand{\BIT}{\begin{itemize}}
\newcommand{\EIT}{\end{itemize}}
\newcommand{\BDI}{\begin{description}}
\newcommand{\EDI}{\end{description}}
\newcommand{\BRM}{\begin{remarks}}
\newcommand{\ERM}{\end{remarks}}

\newcommand{\BTH}{\begin{theo}}
\newcommand{\ETH}{\end{theo}}
\newcommand{\BPR}{\begin{sat}}
\newcommand{\EPR}{\end{sat}}
\newcommand{\BC}{\begin{cases}}
\newcommand{\EC}{\end{cases}}
\newcommand{\COM}[1]{}
\newcommand{\BL}{\begin{lem}}
\newcommand{\EL}{\end{lem}}

\begin{document}

\begin{center}

{\Large \bf Tail Approximation for Reinsurance Portfolios of Gaussian-like Risks}

\vskip 0.4 cm

\centerline{\large Julia Farkas$^1$ and  Enkelejd  Hashorva
\footnote{Department of Actuarial Science, Faculty of Business and Economics,
University of Lausanne, UNIL-Dorigny 1015 Lausanne, Switzerland}}
\vskip 0.4 cm
\centerline{\large University of Lausanne}

\today{}
\end{center}

\textbf{Abstract}: We consider two different portfolios of
proportional reinsurance of the same pool of risks. This contribution is
concerned with Gaussian-like risks, which means that for large values the
survival function of such {risks is, up to a multiplier,} the same as that of
a standard Gaussian risk. We establish the tail asymptotic behavior of the
total loss of each of the reinsurance portfolios and determine also the
relation between randomly scaled Gaussian-like portfolios and unscaled ones.
Further we show that jointly two portfolios of Gaussian-like risks exhibit
asymptotic independence and their weak tail dependence coefficient is non-negative.



\section{Introduction}

In numerous insurance and financial situations the same source of risks
impacts simultaneously different portfolios according to individual
deterministic weights associated with those risks. For instance consider two
big reinsurance companies that operate on the international level, and
{therefore happen to reinsure different proportions of the same risks. If the
reinsurance treaty is a proportional one, then the total risk of each company
for the proportional business is given by a linear combination of risks,
arising from each portfolio of the direct insurer taking part in the
reinsurance programme.}\newline Throughout the paper $X_{i},i\leq n$ will be
independent random variables which alternatively are referred to as risks,
reflecting our interest on insurance and finance applications. In a
probabilistic setting the total loss amount of each reinsurance company can
be modeled by $Q_{n}$ and $W_{n}$, respectively with
\[
Q_{n}=\sum_{i=1}^{n}\lambda_{i}X_{i},\quad W_{n}=\sum_{i=1}^{n}\lambda
_{i}^{\ast}X_{i},
\]
where $X_{i}$ is the financial loss amount claimed from the $i$th direct
insurer, and $\lambda_{i},\lambda_{i}^{\ast}$ are the proportionality factors
of the risks being shared.

In a financial context, as for instance in that considered by Geluk et al.\ (2007), both $Q_{n}$ and $W_{n}$ model {two portfolios with financial returns, where the risks} $X_{1},\ldots,X_{n}$ {are the individual asset returns or
risk factors and $\lambda_{i},\lambda_{i}^{\ast},i\leq n$ are the asset
weights. Typically, the asset weights are assumed to sum to }${1}$.

In concrete insurance and finance applications the distribution function of
financial risks is not known. Essentially, this is not a major drawback, since
often of interest is the quantification of the probability of large
catastrophic risks, especially from the side of the reinsurer. In
applications, departure from a Gaussian model is possible, however for
inference a model with "Gaussian-like" features is of course
preferable.\newline The main purpose of this article is to explore
Gaussian-like risks i.e., risks that are similar to Gaussian ones in terms of
the probability of producing large values. Specifically, we shall assume that
for {any} risk $X_{i},i\leq n$
\begin{equation}
P\left(  X_{i}>u\right)  \sim{\mathcal{L}}_{i}(u)u^{\alpha_{i}}\exp
(-u^{2}/2),\quad u\rightarrow\infty, \label{1}%
\end{equation}
where ${\mathcal{L}}_{i}(\cdot),i\leq n$ are slowly varying functions at
infinity i.e., for any $t$ positive $\lim_{u\rightarrow\infty}{\mathcal{L}%
}_{i}(tu)/{\mathcal{L}}_{i}(u)=1$. In other words
\begin{equation}
P\left(  X_{i}>u\right)  \sim\sqrt{2\pi}{\mathcal{L}}_{i}(u)u^{\alpha_{i}%
+1}\Psi(u),\quad u\rightarrow\infty, \label{1:2}%
\end{equation}
with $\Psi$ the survival function of a $N(0,1)$ random variable; throughout
this paper $f_{1}(x)\sim f_{2}(x)$ means asymptotic equivalence i.e.,
$f_{1}(x)/f_{2}(x)\rightarrow1$ as $x\rightarrow\infty$.

Clearly, if $\alpha_{i}=-1$ and ${\mathcal{L}}_{i}(u)\to(2\pi)^{-1/2}$ as
$u\to\infty$, then $X_{i}$ is tail equivalent to a $N(0,1)$ random variable.
However, in general \rH{a} Gaussian-like risk differs very strongly from \rH{a} Gaussian one
since $\alpha_{i}$ can take large negative values. It is therefore interesting
to investigate the individual behavior {of} each portfolio consisting of
Gaussian-like risks in terms of the probabilities of observing large losses.
We shall investigate first the asymptotic {behavior} of $P\left(
Q_{n}>u\right)  $ for $u\rightarrow\infty$. In view of Lemma 8.6 in Piterbarg
(1996)
\begin{equation}
P\left(  X_{1}+X_{2}>u\right)  \sim{\mathcal{L}}_{1}{\mathcal{L}}_{2}\sqrt
{\pi}\frac{u^{\alpha_{1}+\alpha_{2}+1}}{2^{\alpha_{1}+\alpha_{2}+1}}%
\exp(-u^{2}/4),\quad u\rightarrow\infty\label{PIT}%
\end{equation}
holds for any two Gaussian-like risks $X_{1},X_{2}$ satisfying \eqref{1} with
${\mathcal{L}}_{i}(u)\equiv{\mathcal{L}}_{i}>0,\forall u>0,i=1,2$, which
implies the tail asymptotic behavior of $Q_{2}$ for the case $\lambda
_{1}=\lambda_{2}>0$.\newline However, if $\lambda_{1}\not =\lambda_{2}$ and,
more generally, if risks obey \eqref{1}, then the tail asymptotics of $Q_{2}$
cannot be established by simply using \eqref{PIT}. \rH{The risks}
obeying \eqref{1} do not belong to the class of subexponential distributions
(see Embrechts et al.\ (1997) or Foss et al.\ (2011) for the properties of
this class). In fact {those} risks belong to the class of superexponential
distributions, see Rootz\'{e}n (1986,1987), Kl\"{u}ppelberg and Lindner
(2005), or Geluk et al. (2007) for more details.\newline
We note in passing
that if $X_{1},X_{2}$ are independent $N(0,1)$ random variables, then
$X_{1}+X_{2}$ is a $N(0,2)$ random variable, so \eqref{PIT} follows easily.
Therefore, the appearance of $\exp(-u^{2}/4)$ in the general case in
\eqref{PIT} is intuitively expected since we deal with "Gaussian-like" risks.\\
As it will be discussed below, special Gaussian-like risks \bE{relate} to the random scaling of Gaussian risks.\\
Indeed,  the random scaling is a common phenomena in various insurance models which incorporate inflation or deflation.
In our framework, the random scaling of  $X_i$'s will be modelled by non-negative random variables $S_i,i\le n$ being independent of $X_i,i\le n$. Under certain restrictions, it follows that the randomly scaled risk  $S_i X_i$ is a Gaussian-like one, if $X_i$ is a Gaussian like risk.
This closure property together with the Gaussianity of $X_i$'s are crucial \bE{for} extending \eqref{PIT} to
Gaussian-like risks obeying \eqref{1}. Further, the random scaling technique
utilized in the proof of the main result leads to the derivation of the tail
asymptotic behavior of $Q_{n}$ if each risk $S_{i}$ is bounded, and its
survival function is regularly varying at its upper endpoint, see \eqref{U}
below.\newline Our new result allows us to calculate the weak tail dependence
coefficient $\bar{\chi}(Q_{n},W_{n})$. This measure
of asymptotic independence introduced in Coles et al.\ (1999) is important for
modelling of joint extremes.

The organization of the rest of the paper: we continue below with the
{formulation of the} main results. Section 3 presents two applications. The
first one establishes the asymptotic independence of both portfolios $Q_{n}$
and $W_{n}$, whereas the second one derives the weak tail dependence
coefficient $\bar{\chi}(Q_{n},W_{n})$. \bE{All the proofs are} relegated to Section 4.

\section{Main Results}

In the following $X_{i}$'s are independent (but not identically distributed)
risks with distribution functions $F_{i},i\leq n$ and $S_{i},i\leq n$ are
independent non-negative risks with distribution function $G_{i},i\leq n$. We
shall write for short $X_{i}\sim F_{i},S_{i}\sim G_{i},i\leq n$.  Further, we
shall assume that $X_{1},\ldots,X_{n},S_{1},\ldots,S_{n}$ are mutually
independent. In the special case $X_{i}\sim N(0,1),i\leq n$ and $\lambda
_{i},i\leq n$ are given constants
\[
Q_{n}^{\ast}:=\sum_{i=1}^{n}\lambda_{i}\sqrt{S_{i}}X_{i}\overset{d}{=}%
X_{1}\sqrt{\sum_{i=1}^{n}\lambda_{i}^{2}S_{i}}=:X_{1}\sqrt{V_{n}},
\]
where $\overset{d}{=}$ means equality of distribution functions.\\
For practical applications due to the time-value considerations of money
random scaling of $X_{i}^{\prime}$s by $S_{i}$'s is natural. If as above the
$X_{i}^{\prime}$s are normally distributed, then instead of considering the
tail asymptotics of $Q_{n}$ we can investigate that of $V_{n}$, and $X_{1}$ separately
and then determine the tail asymptotics of $Q_{n}^{*}$. Indeed, by the fact that $X_{1}$ and $V_{n}$ are independent, and
the tail asymptotics of $X_{1}$ is known, in view of Hashorva et al.\ (2010),
the tail asymptotic behavior of the portfolio of risks modeled by $Q_{n}%
^{\ast}$ follows under certain assumptions on $V_{n}$ which are satisfied if
the $G_{i}$'s have a finite upper {endpoint} $\omega_{i}:=\sup(x:G_{i}(x)<1)$
and if $1-G_{i}$ is regularly varying at $\omega_{i}$. More specifically, we
shall assume that $\omega_{i}=1,i\leq n$ and
\begin{equation}
\lim_{x\rightarrow\infty}\frac{P\left(  S_{i}>1-t/x\right)  }{P(S_{i}%
>1-1/x)}=t^{\gamma_{i}},\quad\forall t>0 \label{U}%
\end{equation}
for each $i\leq n$ with some index $\gamma_{i}\in\lbrack0,\infty)$.

Our result on the tail behavior of $V_{n}$ is surprising in that it links the
tail asymptotic behavior of the aggregated risk with that of the products of
the risks. The arithmetic-geometric mean inequality implies that
\[
V_{n}\geq\prod_{i=1}^{n}S_{i}^{\lambda_{i}}=:V_{n}^{\ast},\quad1\leq i\leq n,
\]
provided that $\sum_{i=1}^{n}\lambda_{i}=1$. Our first result below shows the surprising fact that $V_{n}$ and $V_{n}%
^{\ast}$ have the same tail asymptotic behavior.

\begin{theo}
\label{ProdSum} Let $S_{i}\sim G_{i},i\leq n$ be independent non-negative
random variables satisfying \eqref{U}. Then for any $\lambda_{i}>0,i\leq n$
such that $\sum_{i=1}^{n}\lambda_{i}=1$
\[
P\Bigl (\sum_{i=1}^{n}\lambda_{i}S_{i}>u\Bigr )\sim P\Bigl
(\prod_{i=1}^{n}S_{i}^{\lambda_{i}}>u\Bigr )\sim\frac{\prod_{i=1}^{n}%
\lambda_{i}^{-\gamma_{i}}\Gamma(\gamma_{i}+1)P(S_{i}>u)}{\Gamma(\sum_{i=1}%
^{n}\gamma_{i}+1)},\quad u\uparrow1
\]
holds, where $\Gamma(\cdot)$ is the Euler Gamma function.
\end{theo}

Next, we show how by using this theorem, we can reduce the proof of the following
Theorem \ref{kr} in an important particular case to random scaling of a
portfolio of independent standard Gaussian variables. Under the assumptions of
{Theorem \ref{ProdSum}} we know the asymptotic behavior of {$\widetilde
{V_{n}}= \sum_{i=1}^{n} \lambda_{i}^{2} S_{i}/\lambda^{2}, \lambda^{2}=
\sum_{i=1}^{n} \lambda_{i}^{2}$}, and in particular we find that
\[
P \left( \sqrt{ {\widetilde{V_{n}}}}>1-1/u \right) \sim P( {\widetilde{V_{n}}%
}>1-2/u),\quad u \rightarrow\infty.
\]
The distribution function of $X_{1}$ is in the Gumbel max-domain of attraction (MDA)
 with scaling function $w(x)=x$. We recall that a random variable $Z$ with
$P(Z>u)<1,\forall u>0$ is in the Gumbel MDA with some
positive scaling function $w(\cdot)$ if
\[
P(Z>u+s/w(u))\sim\exp(-s)P(Z>u),\quad u\rightarrow\infty
\]
holds for any $s\geq0$, see e.g., Embrechts et al.\ (1997). Applying Theorem
3.1 of Hashorva et al.\ (2010) (see also Hashorva (2012)) we obtain thus
\begin{align}\label{VN}
P(Z_{n}>u)  &  \sim P\Bigl (X_{1}\sqrt{ {\widetilde{V_{n}}}}> {u/\lambda
}\Bigr )  \sim\frac{1}{2}P\Bigl (\lvert X_{1}\rvert\sqrt{ {\widetilde{V_{n}}}}>
{u/\lambda}\Bigr )\nonumber\\
&  \sim\frac{1}{2}\Gamma(\sum_{i=1}^{n}\gamma_{i}+1)P\Bigl (\sqrt{
{\widetilde{V_{n}}}}>1-\lambda^{2}/u^{2}\Bigr )P(\lvert X_{1}\rvert
>u/\lambda)\nonumber\\
&  \sim\Gamma(\sum_{i=1}^{n}\gamma_{i}+1)P( {\widetilde{V_{n}}}>1-2(\lambda
/u)^{2})\Psi(u/\lambda),\quad u\rightarrow\infty.
\end{align}
Since $\lim_{u\rightarrow\infty}P( {\widetilde{V_{n}}}>1-1/u^{2})=0$ the
above works only for $\alpha_{i}+1<0,i\leq n$. Thus the above chain of
asymptotic relations leads us to the main result of this paper in the
particular case $\alpha_{i}<-1,$ that is for the distributions possessing
(\ref{1}) with tails lighter than Gaussian. Next \rH{we state our}
main result for all values $\alpha_{i}\in\mathbb{R}$, $i\leq n$.

\begin{theo}
\label{kr} If $X_{i},i\leq n$ are independent Gaussian-like risks satisfying
\eqref{1} for some $\alpha_{i}\in\mathbb{R}$, $i\leq n,$ then for any set of
deterministic weights $\lambda_{i}>0,$ $i\leq n$ we have
\begin{equation}
P(Q_{n}>u)\sim\frac{(\sqrt{2\pi})^{n-1}\prod_{j=1}^{n}\left[  \lambda
_{j}^{\alpha_{j}+1}\mathcal{L}_{j}(u)\right]  u^{\alpha+n-1}}{\lambda
^{2\alpha+2n-1}}\exp\left(  -\frac{u^{2}}{2\lambda^{2}}\right)  \text{ }
\label{last}%
\end{equation}
as $u\rightarrow\infty,$ where $\lambda^{2}=\sum_{i=1}^{n}\lambda_{i}^{2},$
$\alpha=\sum_{i=1}^{n}\alpha_{i}.$
\end{theo}

\textbf{Remarks}: a) In {Theorem \ref{kr}} we do not put any assumption on the
lower asymptotic tail behavior of the risks. {In the Gaussian mean-zero case
such risks are symmetric about 0. If in Theorem \ref{kr} we assume that the
Gaussian-like risks are symmetric about zero}, then (\ref{last}) can be easily
adapted to the case that $\lambda_{i}\in\mathbb{R},i\leq n$. \\
b) If
${\mathcal{L}}_{i}(\cdot)$ is a constant function, then as mentioned in the
Introduction the risk $X_{i}$ belongs to the class of superexponential
distributions. The tail asymptotics of the convolution of identically
distributed and independent superexponential risks is established in
Rootz\'{e}n (1987) and for more general risks in Kl\"{u}ppelberg and Lindner
(2005). In the aforementioned papers the results are derived under several
constraints on the probability density function of the risks, which we do not
impose here. Hence both our results and Lemma 8.6 of Piterbarg (1996) do not
follow from Rootz\'{e}n (1987) or Kl\"{u}ppelberg and Lindner (2005). \newline
d) In view of {Theorem \ref{ProdSum}} and {Theorem \ref{kr}} the total loss of
randomly scaled risks modeled by $Q_{n}^{\ast}$ is a Gaussian-like risk if the
$S_{i}$'s and $X_{i}^{\prime}$s are independent and satisfy the assumptions of
{Theorem \ref{ProdSum}} and {Theorem \ref{kr}}, respectively.

\bigskip The proof of Theorem \ref{kr} in the general case is based on the
following generalization of Lemma 8.6 of Piterbarg (1996) to random variables
obeying (\ref{1}).

\begin{lem}
\label{8.6} If $X_{i},i=1,2$ are two independent random variables such that
\begin{align}
\label{e86}P\left(  X_{i}> u\right)  {\sim} \mathcal{L}_{i}(u)u^{\alpha_{i}%
}\exp\left(  -\frac{u^{2}}{2p_{i}^{2}}\right)  , \quad u\rightarrow
\infty,\ \ i=1,2
\end{align}
for some $\alpha_{i}\in\mathbb{R}$, $p_{i}>0$ and $\mathcal{L}_{i}(\cdot)$ are
slowly varying functions at infinity, then as $u\rightarrow\infty$
\begin{align}
\label{eex}
P\left(  X_{1}+X_{2}> u\right)  {\sim}\frac{\sqrt{2\pi}p_{1}^{2\alpha_{1}%
+1}p_{2}^{2\alpha_{2}+1}\mathcal{L}_{1}(u)\mathcal{L}_{2}(u)u^{\alpha
_{1}+\alpha_{2}+1}}{p^{2\alpha_{1}+2\alpha_{2}+3}}\exp\left(  -\frac{u^{2}%
}{2p^{2}}\right)  ,
\end{align}
with $p=\sqrt{p_{1}^{2}+ p_{2}^{2}}$.
\end{lem}

\textbf{Example 1}. Consider $X_{1},X_{2}$ two independent Gaussian-like risks
which satisfy \eqref{1}. Applying \eqref{eex} with $p_{1}=p_{2}=1,
p=\sqrt{2}$ we obtain
\begin{align*}
P(X_{1}+X_{2}>u)  &  {\sim} \sqrt{ \pi} \mathcal{L}_{1}(u) \mathcal{L}_{2}(u)
(u/2)^{\alpha_{1}+ \alpha_{2}+1} \exp(-u^{2}/4)
\end{align*}
as $u\to\infty$, which implies \eqref{PIT}.
In particular, when $X_{1},X_{2}$ are independent $N(0,1)$ random variables,
then $X_{1}+X_{2}$ is a $N(0,2)$ random variable, and therefore its tail
asymptotics is given by
\[
P\left(  X_{1}+X_{2}> u\right)  = P\left(  X_{1}> u/\sqrt{2}\right)  \sim
\frac{ 1}{\sqrt{ \pi} u} \exp(- u^{2}/4), \quad u\to\infty,
\]
which follows also from \eqref{e86} when  $\mathcal{L}_{1}=\mathcal{L}_{2}=(2
\pi)^{-1/2}$ and $\alpha_{i}=-1,i=1,2$.

\section{Applications}

A bivariate Gaussian random vector $(X,Y)$ with $N(0,1)$ marginals is
specified completely by the correlation coefficient $\rho$. Although $\rho<1$
can be very close to $1$, still $X$ and $Y$ are asymptotically independent in
the sense that
\begin{equation}
\lim_{u\rightarrow\infty}\frac{P(X>u,Y>u)}{P(X>u)}=0. \label{pje}%
\end{equation}
Asymptotic independence is a nice property, closely related to joint
asymptotic behavior of componentwise sample maxima (e.g., Resnick (1987)). For
bivariate Gaussian samples the componentwise maxima are (using (\ref{pje}))
asymptotically independent. The asymptotic independence is a crucial property
for the calculation of many indices related to extreme value statistics,
finance and insurance applications. In our first application we show that
\rH{the losses modeled by} $Q_{n}$ and $W_{n}$ are asymptotically independent.\newline If
$Q_{n}$ and $W_{n}$ have distribution functions $H$ and $H_{\ast}$,
respectively, then the asymptotic independence of $Q_{n}$ and $W_{n}$ means
that
\[
\chi_{u}(Q_{n},W_{n}):=\frac{P(Q_{n}>t_{u},W_{n}>t_{u}^{\ast})}{P(Q_{n}%
>t_{u})}\rightarrow0,\quad u\rightarrow\infty,
\]
with $t_{u}:=H^{-1}(1-1/u),t_{u}^{\ast}:=H_{\ast}^{-1}(1-1/u),$ and
$H^{-1},H_{\ast}^{-1}$ are the generalized inverses of $H$ and $H_{\ast}$,
respectively. \COM{ Clearly, the relation $\lim_{u\rightarrow\infty}\chi
_{u}(Q_{n}/a,W_{n}/b)=0$ for some positive $a$ and $b$ is equivalent to
$\lim_{u\rightarrow\infty}\chi_{u}(Q_{n},W_{n})=0$. Choose thus \[
a=\sqrt{\sum_i=1^n\lambda_i^2},\quad b=\sqrt{\sum_i=1^n(\lambda_i^\ast)^2.}\]
Now it is convenient instead of assuming
\begin{equation}
\sum_i=1^n\lambda_i=\sum_i=1^n\lambda_i^\ast=1, \label{basic}%
\end{equation}
assume that $\lambda_{i},\lambda_{i}^{\ast},i\leq n$ are positive and satisfy
\begin{equation}
\sum_i=1^n\lambda_i^2=\sum_i=1^n(\lambda_i^\ast)^2=1. \label{both}%
\end{equation}
This gives no restriction to our analysis, since one can easily turn back from
(\ref{both}) to (\ref{basic}), see the remark below Theorem \ref{KK}. The
assumption (\ref{both}) is reasonable since when all $X_{i}$'s are $N(0,1)$
random variables, then $Q_{n}$ and $W_{n}$ are also $N(0,1)$ distributed. Note
that (\ref{both}) implies \[ \varrho:=\sum_i=1^n\lambda_i\lambda_i^\ast
\in\lbrack0,1]. \] } In view of {Theorem \ref{kr}} both $Q_{n}$ and $W_{n}$
have distribution function in the Gumbel MDA. Utilizing the formula for the
norming constants of Weibull-like distributions given on p. 317 of Mikosch
(2009) and using again Theorem \ref{kr}, it follows that  {with $\lambda=
\sqrt{ \sum_{i=1}^{n} \lambda_{i}^{2}} >0, \lambda^{*}= \sqrt{\sum_{i=1}^{n}
(\lambda_{i}^{*})^{2}}>0$}
\begin{equation}
{\frac{t_{u}}{\lambda} \sim\frac{t_{u}^{\ast}}{\lambda^{*}}},\quad
u\rightarrow\infty.
 \label{HH}%
\end{equation}
{Consequently, since for all $u$ large
\[
\chi_{u}(Q_{n},W_{n})=\frac{P(Q_{n}/\lambda>t_{u}/\lambda,W_{n}/\lambda
>t_{u}^{\ast} /\lambda^{*} )}{P(Q_{n} / \lambda>t_{u}/\lambda)}
\]
we shall assume without loss of generality that} $\lambda_{i},\lambda
_{i}^{\ast},i\leq n$ are positive and satisfy
\begin{equation}
\sum_{i=1}^{n}\lambda_{i}^{2}=\sum_{i=1}^{n}(\lambda_{i}^{\ast})^{2}=1.
\label{both}%
\end{equation}
The assumption (\ref{both}) is reasonable since when all $X_{i}$'s are
$N(0,1)$ random variables, then $Q_{n}$ and $W_{n}$ are also $N(0,1)$
distributed. Note that (\ref{both}) implies
\[
\varrho:=\sum_{i=1}^{n}\lambda_{i}\lambda_{i}^{\ast}\in\lbrack0,1].
\]

Since both portfolios are supposed to be different, we shall assume below that
\begin{equation}
\varrho\in\lbrack0,1). \label{ro}%
\end{equation}
Hence, by {Theorem \ref{kr}}, \eqref{HH} and \eqref{ro} for any $\varepsilon
>0$ we obtain
\begin{align*}
\chi_{u}(Q_{n},W_{n})  &  \leq 
\frac{P(\sum_{i=1}^{n}(\lambda_{i}+\lambda_{i}^{\ast})X_{i}>2t_{u}%
(1+o(1)))}{P(W_{n}>t_{u})}\rightarrow0,\quad u\rightarrow\infty.
\end{align*}
When asymptotic independence holds, as suggested by Coles et al.\ (1999) more
insight on the strength of the joint tail behavior is obtained by calculating
the weak tail dependence coefficient $\bar{\chi}(Q_{n},W_{n}):=\lim
_{u\rightarrow\infty}\bar{\chi}_{u}(Q_{n},W_{n})$ (supposing the limit
exists), where
\[
\bar{\chi}_{u}(Q_{n},W_{n})=\frac{\ln P(Q_{n}>t_{u})+\ln P(W_{n}>t_{u}^{\ast
})}{\ln P(Q_{n}>t_{u},W_{n}>t_{u}^{\ast})}-1.
\]
Borrowing the idea of Piterbarg and Stamatovic (2005) (see also D\c{e}bicki et
al.\ (2010)) we have for all large $u$%
\[
P(Q_{n}>t_{u},W_{n}>t_{u}^{\ast})\leq\inf_{a,b>0,a+b=1}P\Bigl (\sum_{i=1}%
^{n}(a\lambda_{i}+b\lambda_{i}^{\ast})X_{i}>t_{u}(1+o(1))\Bigr ),
\]
which implies using \bE{further} Theorem \ref{kr}
\begin{align}
\limsup_{u\rightarrow\infty}\frac{\ln P(Q_{n}>t_{u},W_{n}>t_{u}^{\ast})}%
{t_{u}^{2}} &  \leq\lim_{u\rightarrow\infty}\frac{\inf_{{a,b\in\mathbb{R}%
,\rH{a+b}=1}}\ln\Bigl(P\Bigl (\sum_{i=1}^{n}(a\lambda_{i}+b\lambda
_{i}^{\ast})X_{i}>t_{u}(1+o(1))\Bigr )\Bigl)}{t_{u}^{2}}\nonumber\\
&  =-\frac{1}{1+\varrho}.\label{PS}%
\end{align}
Consequently,
\[
\limsup_{u\rightarrow\infty}\bar{\chi}_{u}(Q_{n},W_{n})\leq\varrho.
\]
Our last result shows that $\bar{\chi}(Q_{n},W_{n})=\varrho$ for two
Gaussian-like portfolios $Q_{n}$ and $W_{n}$.

\begin{theo}
\label{KK} Under the assumptions of {Theorem \ref{kr}}, if further
\eqref{both} and \eqref{ro} are satisfied, then
\begin{equation}
\bar{\chi}(Q_{n},W_{n})=\varrho\label{rs}%
\end{equation}
holds. Moreover, (\ref{rs}) still holds even if some $\lambda_{i},\lambda
_{i}^{\ast}$ equal zero.
\end{theo}

\textbf{Remark:} If we do not assume (\ref{both}) in {Theorem \ref{KK}}, then
(\ref{rs}) is valid with $\varrho=\frac{\sum
_{i=1}^{n}\lambda_{i}\lambda_{i}^{\ast}}{\sqrt{\sum_{i=1}^{n}\lambda_{i}%
^{2}\sum_{i=1}^{n}(\lambda_{i}^{\ast})^{2}}}$, provided that (\ref{ro})
holds.\newline

\section{Proofs}

\textsc{Proof of Theorem} \ref{ProdSum}. For all $x$ large by the independence
of $S_{1},\ldots,S_{n}$ we may write (set $G_{i,x}(z):=G_{i}(1-z/x),x,z\in
(0,\infty)$)
\begin{align*}
\lefteqn{P(\sum_{i=1}^{n}\lambda_{i}S_{i}>1-1/x)}\\
&  =\int_{0}^{1}P(\lambda_{1}S_{1}>1-1/x-\sum_{i=2}^{n}\lambda_{i}%
y_{i})\,dG_{2}(y_{2})\cdots dG_{n}(y_{n})\\
&  =\int_{0}^{1}P(\lambda_{1}S_{1}>\lambda_{1}-1/x+\sum_{i=2}^{n}\lambda
_{i}-\sum_{i=2}^{n}\lambda_{i}(1-z_{i}/x))\,dG_{2}(z_{2})\cdots dG_{n}%
(z_{n})\\
&  =\int_{0}^{\infty}\cdots\int_{0}^{\infty}P(S_{1}>1-(1-\sum_{i=2}^{n}%
\lambda_{i}z_{i})/(\lambda_{1}x))\,dG_{2,x}(z_{2})\cdots dG_{n,x}(z_{n})\\
&  =\prod_{i=1}^{n}\overline{G}_{i}(1-1/x)\int_{0}^{\infty}\cdots\int
_{0}^{\infty}\frac{P(S_{1}>1-(1-\sum_{i=2}^{n}\lambda_{i}z_{i})/(\lambda
_{1}x))}{\overline{G}_{1}(1-1/x)}\,dG_{2,x}(z_{2})\cdots dG_{n,x}%
(z_{n})/(\prod_{i=2}^{n}\overline{G}_{i}(1-1/x)).
\end{align*}
Assume for simplicity that $n>2$. By the assumption on $G_{i}$, for any
$z_{i}>0,i\leq n$
\[
\frac{\overline{G}_{i,x}(z_{i})}{\overline{G}_{i}(1-1/x)}=\frac{\overline
{G}_{i}(1-z_{i}/x)}{\overline{G}_{i}(1-1/x)}\rightarrow z_{i}^{\gamma_{i}%
},\quad x\rightarrow\infty
\]
and
\[
\lim_{x\rightarrow\infty}\frac{P(S_{1}>1-(1-\sum_{i=2}^{n}\lambda_{i}%
z_{i})/(\lambda_{1}x))}{\overline{G}_{1}(1-1/x)}=\Bigl(\max\Bigl(0,1-\sum
_{i=2}^{n}\lambda_{i}z_{i}\Bigr)\Bigr)^{\gamma_{1}},
\]
which implies as $x\rightarrow\infty$
\begin{align*}
\lefteqn{P(\sum_{i=1}^{n}\lambda_{i}S_{i}>1-1/x)}\\
&  \sim\prod_{i=1}^{n}\overline{G}_{i}(1-1/x)\lambda_{1}^{-\gamma_{1}}%
\prod_{i=2}^{n}\gamma_{i}\int_{0}^{\infty}\cdots\int_{0}^{\infty}%
\Bigl(\max\Bigl(0,1-\sum_{i=2}^{n}\lambda_{i}z_{i}\Bigr)\Bigr)^{\gamma_{1}%
}\prod_{i=2}^{n}z_{i}^{\gamma_{i}-1}\,dz_{2}\cdots dz_{n}\\
&  =\frac{1}{\Gamma(\sum_{i=1}^{n}\gamma_{i}+1)}\prod_{i=1}^{n}\Bigl(\lambda
_{i}^{-\gamma_{i}}\Gamma(\gamma_{i}+1)\overline{G}_{i}(1-1/x)\Bigr).
\end{align*}
Applying Theorem 3.1 in Hashorva et al.\ (2010) we obtain as $x\rightarrow
\infty$
\begin{align*}
\Gamma(\sum_{i=1}^{n}\gamma_{i}+1)P\Bigl (\prod_{i=1}^{n}S_{i}^{\lambda_{i}%
}>1-1/x\Bigr )  &  \sim\prod_{i=1}^{n}\Bigl(P(S_{i}^{\lambda_{i}}%
>1-1/x)\Gamma(\gamma_{i}+1)\Bigr)\\
&  =\prod_{i=1}^{n}\Bigl(P(S_{i}>1-1/(\lambda_{i}x))\Gamma(\gamma
_{i}+1)\Bigr),
\end{align*}
hence the proof is complete.

\bigskip

\textsc{Proof of Theorem} \ref{kr}. In light of Lemma \ref{8.6} for random
variables $\lambda_{i}X_{i},$ $i=1,\ldots,n$ we have $p_{i}^{2}=\lambda
_{i}^{2}$ and correspondingly scaled $\mathcal{L}_{i}$'s. Thus for $n=2$ we have proven that%
\[
P\left(  Q_{2}>u\right)  =\frac{\sqrt{2\pi}\lambda_{1}^{\alpha_{1}+1}%
\lambda_{2}^{\alpha_{2}+1}\mathcal{L}_{1}(u)\mathcal{L}_{2}(u)u^{\alpha
_{1}+\alpha_{2}+1}}{\theta_{2}^{2\alpha_{1}+2\alpha_{2}+3}}\exp\left(
-\frac{u^{2}}{2\theta_{2}}\right)  (1+o(1))
\]
as {$u\rightarrow\infty,$} with $\theta_{2}=\lambda_{1}^{2}+\lambda_{2}^{2}.$
Now we proceed by induction assuming that
\[
P\left(  Q_{k}>u\right)  \sim\frac{(\sqrt{2\pi})^{k-1}\prod_{j=1}^{k}\left[
\lambda_{j}^{\alpha_{j}+1}\mathcal{L}_{j}(u)\right]  u^{\alpha_{1}%
+\cdots+\alpha_{k}+k-1}}{\theta_{k}^{2\alpha_{1}+\cdots+2\alpha_{k}+2k-1}}%
\exp\left(  -\frac{u^{2}}{2\theta_{k}}\right)  ,\text{ }%
\]
with $\theta_{k}=\lambda_{1}^{2}+\cdots+\lambda_{k}^{2}$ and $k>2$.
Considering that $Q_{k+1}=Q_{k}+\lambda_{k+1}X_{k+1}$ with $Q_{k}$ being
independent of $X_{k+1}$, the claim follows by a direct application of Lemma
\ref{8.6}. \hfill$\Box$

\textsc{Proof of Lemma} \ref{8.6}. Let $F_{i},i=1,2$ denote the distribution
functions of \rH{$X_1$ and $X_2$}, respectively. Suppose without loss of generality that $p_{1}^{2}+p_{2}^{2}=1$
and $p_{1}\leq p_{2}$. Then for $c=1.1$ and any $\varepsilon>0$ we have
\begin{align}
P({X_{1}+X_{2}} &  {>x,X_{2}\leq(1-cp_{1})x)}\leq P({X_{1}>cp_{1}%
x)}=O(x^{\alpha_{1}+\varepsilon}\exp(-c^{2}x^{2}/2)),\nonumber\\
P({X_{1}+X_{2}} &  {>x,X_{2}>cp_{2}x)}\leq P({X_{2}>cp_{2}x)}=O(x^{\alpha
_{2}+\varepsilon}\exp(-c^{2}x^{2}/2))\label{8.22}%
\end{align}
as $x\rightarrow\infty.$ Note that $0<a:=1-cp_{1}<b:=cp_{2};$ the first
inequality follows from $p_{1}\leq1/\sqrt{2}$ and the second one follows from
$p_{1}+p_{2}>1.$ Let us focus on the asymptotic behavior of the integral%
\[
I_{x}=\int_{ax}^{bx}\overline{F_{1}}(x-y))dF_{2}(y),\quad\overline{F_{i}%
}=1-F_{i}.
\]
Pick small $h>0$ and denote $h_{k}=kh/x,$ $\Delta_{k}=[h_{k},h_{k+1}),$ where
$k$ is a positive integer. Then, bounding $F_{1}$ on intervals $[x-h_{k}%
,x-h_{k-1}]$ by its maximum and minimum values, respectively and then
integrating in $y$ we have
\[
\sum_{k:\ \Delta_{k}\subset\lbrack ax,bx]}\overline{F_{1}}(x-h_{k-1}%
)(F_{2}(h_{k})-F_{2}(h_{k-1}))\leq I_{x}\leq\sum_{k:\ \Delta_{k}\cap\lbrack
ax,bx]\neq\emptyset}\overline{F_{1}}(x-h_{k})(F_{2}(h_{k})-F_{2}(h_{k-1})).
\]
Observe that there exist two positive functions $A_{1},A_{2}$ decreasing to
zero as $x\rightarrow\infty$ such that for $i=1,2$
\[
\mathcal{L}_{i}(x)x^{\alpha_{i}}\exp(-x^{2}/2p_{i}^{2})(1-A_{i}(x))\leq
\overline{F_{i}}(x)\leq\mathcal{L}_{i}(x)x^{\alpha_{i}}\exp(-x^{2}/2p_{i}%
^{2})(1+A_{i}(x)),\quad\forall x>0.
\]
Similarly, there exist two positive {functions} $B_{1},B_{2}$ decreasing to
zero as $x\rightarrow\infty$ such that
\[
1-B_{i}(x)\leq\inf_{y\in\lbrack a,b]}\frac{\mathcal{L}_{i}(xy)}{\mathcal{L}%
_{i}(x)}\leq\sup_{y\in\lbrack a,b]}\frac{\mathcal{L}_{i}(xy)}{\mathcal{L}%
_{i}(x)}\leq1+B_{i}(x),\quad i=1,2.
\]
Since $x^{\alpha}e^{-x^{2}/2q^{2}},q>0$ decreases for all
sufficiently large $x$ denoting
\[
\gamma_{2}(x)=A_{2}(ax)+B_{2}(x)+A_{2}(ax)B_{2}(x),\ \ r(x)=\frac{e^{x}-1}{x}
\]
we obtain 
\begin{align*}
F_{2}(h_{k}) &  -F_{2}(h_{k-1})
  \leq\mathcal{L}_{2}(x)\Biggl[r(bh/p_{2}^{2})h_{k}^{\alpha_{2}}e^{-h_{k}%
^{2}/2p_{2}^{2}}\frac{kh^{2}}{p_{2}^{2}x^{2}}+2\text{ }\gamma_{2}%
(x)(x)h_{k-1}^{\alpha_{2}}e^{-h_{k-1}^{2}/2p_{2}^{2}}\Biggr].
\end{align*}
In order to derive an estimation from below, note that for sufficiently large $x$
\begin{align*}
\frac{(2k-1)h^{2}}{2x^{2}p_{2}^{2}}+\log\left(  \frac{k-1}{k}\right)
^{\alpha_{2}}
\ge &  \frac{kh^{2}}{p_{2}^{2}x^{2}}(1-C/x^{2})
\end{align*}
for some $C>0$ which does not depend on $h$ for all sufficiently small $h.$ Therefore
\[
F_{2}(h_{k})-F_{2}(h_{k-1})\geq\mathcal{L}_{2}(x)\Biggl[h_{k}^{\alpha_{2}%
}e^{-h_{k}^{2}/2p_{2}^{2}}\frac{kh^{2}}{p_{2}^{2}x^{2}}r\left(  \frac
{kh^{2}(1-C/x^{2})}{p_{2}^{2}x^{2}}\right)  (1-C/x^{2})-2\gamma_{2}%
(x)(x)h_{k-1}^{\alpha_{2}}e^{-h_{k-1}^{2}/2p_{2}^{2}}\Biggr].
\]

Thus we have
\begin{align}
I_{x}  &  \leq\mathcal{L}_{1}(x)\mathcal{L}_{2}(x)\Biggl[(1+A_{1}%
(x/4))(1+B_{1}(x/4))(1+\gamma_{2}(x))r\left(  hb/p_{2}^{2}\right) \nonumber\\
&  \times\sum_{k:\ \Delta_{k}\cap\lbrack ax,bx]\neq\emptyset}(x-h_{k}%
)^{\alpha_{1}}\exp(-(x-h_{k})^{2}/2p_{1}^{2})h_{k}^{\alpha_{2}}e^{-h_{k}%
^{2}/2p_{2}^{2}}\frac{kh^{2}p_{2}^{2}}{x^{2}}\nonumber\\
&  +2(1+A_{1}(x/4))\gamma_{2}(x)(1+B_{1}(x))\nonumber\\
&  \times\sum_{k:\ \Delta_{k}\cap\lbrack ax,bx]\neq\emptyset}(x-h_{k}%
)^{\alpha_{1}}\exp(-(x-h_{k})^{2}/2p_{1}^{2})h_{k-1}^{\alpha_{2}}%
e^{-h_{k-1}^{2}/2p_{2}^{2}}\Biggr] \label{8.25}%
\end{align}
and%
\begin{align}
I_{x}  &  \geq\mathcal{L}_{1}(x)\mathcal{L}_{2}(x)\Biggl[(1-A_{1}%
(x/4))(1-B_{1}(x/4))(1-\gamma_{2}(x))r\left(  ha(1-Cx^{-2})/p_{2}^{2}\right)
\nonumber\\
&  \times\sum_{k:\ \Delta_{k}\subset\lbrack ax,bx]}(x-h_{k})^{\alpha_{1}}%
\exp(-(x-h_{k})^{2}/2p_{1}^{2})h_{k}^{\alpha_{2}}e^{-h_{k}^{2}/2p_{2}^{2}%
}\frac{kh^{2}p_{2}^{2}}{x^{2}}\nonumber\\
&  -2(1+A_{1}(x/4))\gamma_{2}(x/4)(1+B_{1}(x/4))(1-Cx^{-2})\nonumber\\
&  \times\sum_{k:\ \Delta_{k}\cap\lbrack ax,bx]\neq\emptyset}(x-h_{k}%
)^{\alpha_{1}}\exp(-(x-h_{k})^{2}/2p_{1}^{2})h_{k-1}^{\alpha_{2}}%
e^{-h_{k-1}^{2}/2p_{2}^{2}}\Biggr]. \label{8.26}%
\end{align}
The first sums in (\ref{8.25}) and (\ref{8.26}) differ from each other by two
summands, so it is sufficient to estimate one of them. Then the first sum in
the right-hand side of (\ref{8.25}) is equal to (set $h_{k}^{\prime}%
=h_{k}/x=hk/x^{2}$)
\begin{align*}
I_{x}^{\prime}  &  :=p_{2}^{-2}x^{\alpha_{1}+\alpha_{2}+2}\sum_{k:\ \Delta
_{k}/x\cap\lbrack a,b]\neq\emptyset}(1-h_{k}^{\prime})^{\alpha_{1}}%
(h_{k}^{\prime})^{\alpha_{2}}\exp\left(  -\frac{x^{2}(1-h_{k}^{\prime2}%
)}{2p_{1^{2}}}-\frac{x^{2}h_{k}^{2}}{2p_{2}^{2}}\right)  h_{k}^{\prime}%
\frac{h}{x^{2}}\\
&  \leq
p_{2}^{-2}x^{\alpha_{1}+\alpha_{2}+2}\int_{a}^{b}(1-t+h/x^{2})^{\alpha
_{1}}t^{\alpha_{2}+1}\exp\left(  -\frac{x^{2}}{2}\left(  \frac{(1-t)^{2}%
}{p_{1}^{2}}+\frac{t^{2}}{p_{2}^{2}}\right)  \right)  dt,
\end{align*}
where we used the monotonicity of the involved functions. In order to obtain a lower
bound for the first sum in (\ref{8.26}) replace  $(1-t+h/x^{2}%
)^{\alpha_{1}}t^{\alpha_{2}+1}$ by $(1-t)^{\alpha_{1}}(t-h/x^{2})^{\alpha
_{2}+1}.$ Next,
Theorem 1.3 {in} Fedoryuk (1987) \bE{yields}
\[
I_{x}^{\prime}=\sqrt{2\pi}p_{1}^{2a_{1}+1}p_{2}^{2a_{2}+1}x^{\alpha_{1}%
+\alpha_{2}+1}e^{-x^{2}/2}(1+O(x^{-2})),\quad x\rightarrow\infty.
\]
We investigate below the second sums $I_{x}^{\prime\prime}$ and $J_{x}^{\prime\prime
}$ on the right-hand {side} of (\ref{8.25}) and (\ref{8.26}), respectively.
For any $k,$ the $k$th summands in those sums are equal to the $k$th summands
in the first sums multiplied by $x^{2}/(kh)^{2},$ which is not greater than
$b/h.$ Thus we obtain dividing right- and left- parts of (\ref{8.25}) and
(\ref{8.26}) by%
\[
D(x)=\sqrt{2\pi}p_{1}^{2a_{1}+1}p_{2}^{2a_{2}+1}\mathcal{L}_{1}(x)\mathcal{L}%
_{2}(x)x^{\alpha_{1}+\alpha_{2}+1}e^{-x^{2}/2}
\]
and letting $x\rightarrow\infty,$ that%
\[
r\left(  ha/p_{2}^{2}\right)  \leq\liminf_{x\rightarrow\infty}\frac{I_{x}%
}{D(x)}\leq\limsup_{x\rightarrow\infty}\frac{I_{x}}{D(x)}\leq r\left(
hb/p_{2}^{2}\right)  ,
\]
which by definition of $r(x)$ and the arbitrary choice of $h$ establishes the
asymptotic behavior of $I_{x}$. Clearly, in view of the fact that
$\mathcal{L}_{i}(x/p)\sim\mathcal{L}_{i}(x),$ $i=1,2$ the proof in the case
that $p_{1}^{2}+p_{2}^{2}=1$ follows from (\ref{8.22}). The general case of
$p_{1},p_{2}$ follows by re-scaling, and thus the proof is complete.
\hfill$\Box$

\textsc{Proof of Theorem } \ref{KK}. In view of (\ref{PS}) we need to estimate
$P\left(  Q_{n}>t_{u},W_{n}>t_{u}^{\ast}\right)  $ from below. \COM{ We shall
construct a set $A$ such that Introduce a partition of $\mathbb{R}^{n},$\[
\mathcal{R}_u:=\{\mathbf{\Delta}_\mathbf{k}=\mathbf{k\Delta}_\mathbf
{0},\mathbf{k\in}\mathbb{Z}^n\},\ \mathbf{\Delta}_\mathbf{0}=[0,\Delta
_1(u))\times...\times\lbrack0,\Delta_n(u)), \] we shall choose the values for
positive $\Delta_{i}(u)$  later on. We have, denoting $\mathbf{X=(}%
X_{1},...,X_{n}),$%
\begin{align*}
P\left(      Q_n>u,W_n>u\right)       &  =\sum_\mathbf{k\in}\mathbb
{Z}^nP\left(  Q_n>u,W_n>u,\mathbf{X\in\Delta}_\mathbf{k}\right)      \\ &
\geq\sum_\mathbf{k:\Delta}_\mathbf{k}\in\{\Sigma_i=1^n\lambda_ix_i\geq
u,\Sigma_i=1^n\lambda_i^\ast x_i\geq u\}P(\mathbf{X\in\Delta}_\mathbf{k})\\ &
\geq P(\mathbf{X\in\Delta}_\mathbf{k}_0),
\end{align*}
where $\mathbf{k}_{0}$ will be chosen in such a way that $\mathbf{\Delta
}_{\mathbf{k}_{0}}\subset\{\Sigma_{i=1}^{n}\lambda_{i}x_{i}\geq u,\Sigma
_{i=1}^{n}\lambda_{i}^{\ast}x_{i}\geq u\},$ and the asymptotic behavior of the
logarithm of the last probability coincides with (\ref{PS}).} We shall
determine optimal $\delta_{i}(u),i\leq n$ {such that}
\begin{align}
\label{determine}P\left(  Q_{n}>t_{u},W_{n}>t_{u}^{\ast}\right)  \geq P\left(
X_{i}>\delta_{i}(u),i=1,\ldots,n\right)  .
\end{align}
In order to realize such a choice, consider the asymptotic behavior of the
integral%
\[
\int_{\{\Sigma_{i=1}^{n}\lambda_{i}x_{i}\geq u,\Sigma_{i=1}^{n}\lambda
_{i}^{\ast}x_{i}\geq u\}}e^{-\frac{1}{2}\left\Vert \mathbf{x}\right\Vert ^{2}%
}d\mathbf{x=}u^{n}\int_{\{\Sigma_{i=1}^{n}\lambda_{i}s_{i}\geq1,\Sigma
_{i=1}^{n}\lambda_{i}^{\ast}s_{i}\geq1\}}e^{-\frac{1}{2}u^{2}\left\Vert
\mathbf{s}\right\Vert ^{2}}d\mathbf{s.}%
\]
In the spirit of the Laplace asymptotic method, we find the minimal value of $\left\Vert
\mathbf{s}\right\Vert ^{2}=\sum_{i=1}^{n}s_{i}^{2}$ on the set $\{\mathbf{s:}%
\Sigma_{i=1}^{n}\lambda_{i}s_{i}\geq1,\Sigma_{i=1}^{n}\lambda_{i}^{\ast}%
s_{i}\geq1\}.$ Since $\lambda_{i},\lambda_{i}^{\ast}$ are all non-negative,
the minimum is attained at the boundary, that is, on the set $\{\mathbf{s:}%
\Sigma_{i=1}^{n}\lambda_{i}s_{i}=1,\Sigma_{i=1}^{n}\lambda_{i}^{\ast}%
s_{i}=1\}.$
 It follows that the point of minimum \rH{has components}
\[
s_{i}=\frac{\lambda_{i}+\lambda_{i}^{\ast}}{1+\rho},\ i=1,\ldots,n.
\]
Consequently, the minimal value of $\left\Vert \mathbf{s}\right\Vert ^{2}$ on the
integrating set equals $2/(1+\rho).$ Setting now (write $z_{u}:=\max
(t_{u},t_{u}^{\ast})$)
\[
\delta_{i}(u)=\frac{\lambda_{i}+\lambda_{i}^{\ast}}{1+\rho}z_{u}
\]
we have that \eqref{determine} holds for any $u>0$ and furthermore, by (1) and
(9)
\begin{align*}
\log P\left(  X_{i}>\delta_{i}(u),i=1,\ldots,n\right)   &  =\sum_{i=1}^{n}\log
P\left(  X_{i}>\delta_{i}(u)\right) \\
&  \sim-z_{u}^{2}\sum_{i=1}^{n}\frac{(\rho_{i}+\rho_{i}^{\ast})^{2}%
}{2(1+\varrho)^{2}},\quad u\rightarrow\infty\\
 & =-\frac{z_{u}^{2}}{1+\rho}%
\end{align*}
and thus the claim follows using \eqref{HH}.\hfill$\Box$

\textbf{Acknowledgments.} We are in debt to both reviewers for numerous
suggestions which improved our manuscript. We would like to thank also Zakhar
Kabluchko for discussions on related results. E. Hashorva kindly acknowledges
partial support by the Swiss National Science Foundation Grants
200021-1401633/1 and 200021-134785. {J. Farkas has been partially supported by the
project RARE -318984  (a Marie Curie IRSES Fellowship within the 7th European Community Framework Programme)}.

\end{document}